\documentclass[10pt,a4paper,oneside]{amsart}
\usepackage[utf8]{inputenc}
\usepackage[english]{babel}
\usepackage{amsmath}
\usepackage{amsfonts}
\usepackage{amssymb}
\usepackage{lmodern}
\usepackage{amsthm}

\usepackage{color}

\title[Vanishing theorems for free boundary submanifolds]{Vanishing theorems for the cohomology groups of free boundary submanifolds}
\author{Marcos P. Cavalcante}
\author{Abraão Mendes}
\author{Feliciano Vitório}
\date{\today}
\address{Instituto de Matemática, Universidade Federal de Alagoas, Maceió, AL, 57072-970, Brazil}
\email{marcos@pos.mat.ufal.br}
\email{abraao.mendes@im.ufal.br}
\email{feliciano@pos.mat.ufal.br}

\newcommand{\R}{\mathbb{R}}
\newcommand{\ds}{\displaystyle}

\newcommand{\p}{\partial}
\newcommand{\Har}{\mathcal{H}}

\newtheorem{theorem}{Theorem}[section]
\newtheorem{lemma}[theorem]{Lemma}
\newtheorem{corollary}[theorem]{Corollary}

\newtheorem{claim}[theorem]{Claim}

\theoremstyle{remark}

\newcommand{\B}{\mathbb{B}}

\begin{document}

\begin{abstract}
In this paper, we prove that there exists a universal constant $C$, depending only on positive integers $n\geq 3$ and $p\le n-1$, such that if $M^n$ is a compact free boundary submanifold of dimension $n$ immersed in the Euclidean unit ball $\mathbb{B}^{n+k}$ whose  size of the traceless second fundamental form is less than $C$, then the $p$th cohomology group of $M^n$ vanishes. Also, employing a different technique, we obtain a rigidity result for compact free boundary surfaces minimally immersed in the unit ball $\mathbb{B}^{2+k}$.
\end{abstract}

\maketitle

\section{Introduction}

In 1968, Simons \cite{Simons} proved that if $M^n$ is a compact $n$-manifold minimally immersed in the unit sphere $\mathbb S^{n+k}$ whose second fundamental form $A$ satisfies $\|A\|^2\leq \frac{nk}{2k-1}$, then either $\|A\|^2=0$ (i.e. $M^n$ is totally geodesic) or $\|A\|^2=\frac{nk}{2k-1}$. Later, Lawson \cite{Lawson} and Chern, do Carmo, and Kobayashi \cite{CdCK} classified all minimal submanifolds in $\mathbb S^{n+k}$ which satisfy $\|A\|^2=\frac{nk}{2k-1}$. Such submanifolds are either the Veronese surface in $\mathbb S^4$ or a family of products of two spheres with appropriate radii, which are currently known  as minimal Clifford tori. In particular, in codimension one, only Clifford tori occur. These important results say that there exists a gap in the space of minimal submanifolds in $\mathbb S^{n+k}$ in terms of the lenght of their second fundamental forms and their dimensions. This kind of behavior has been observed in many other cases as we can see, for instance, in  \cite{AdC}, \cite{LiLi}, \cite{Flaherty}, \cite{HowardWei}, \cite{LawsonSimons}, \cite{Leung},  \cite{Leung2}, \cite{ShiohamaXu}, \cite{Xu} and references therein.

We point out here the important contribution done by Lawson and Simons in \cite{LawsonSimons}, where they proved a topological gap result without making any minimality assumption on the submanifold. The result is the following.

\begin{theorem}[Lawson-Simons] 
Let $M^n$ be a closed $n$-manifold immersed in the Euclidean unit $(n+k)$-sphere with second fundamental form $A$ satisfying
$$
\|A\|^2 <\min\{p(n-p), 2\sqrt{p(n-p)} \},
$$
where $p\leq n-1$ is a positive integer. Then, for any finitely generated Abelian group $G$, 
$$
H_p(M;G) = H_{n-p}(M;G)=0.
$$
In particular, if $\|A\|^2 <\min\{n-1, 2\sqrt{n-1} \}$, then $M^n$ is a homotopy sphere.
\end{theorem}

 
Taking a new perspective, Ambrozio and Nunes \cite{AmbrozioNunes} recently obtained a geometric gap type theorem for free boundary minimal surfaces $M^2$ in the Euclidean unit $3$-ball $\mathbb{B}^3$. They proved that if $\|A\|^2(x)\langle x, N(x)\rangle^2 \leq 2$, where $N(x)$ is the unit normal vector at $x\in M^2$, then $M^2$ is either the equatorial disk or the critical catenoid. 

We recall that a submanifold $M^n$, with nonempty boundary $\partial M$, which is minimally immersed in the unit ball $\mathbb B^{n+k}$ and such that $M\cap\partial\mathbb B^{n+k} = \partial M$ is called \emph{free boundary} if $M^n$ intersects $\partial \mathbb B^{n+k} = \mathbb S^{n+k-1}$ in a right angle along its boundary $\partial M$. Such submanifolds are critical points for the area functional for those variations that keep the boundary of $M^n$ into the boundary of $\mathbb B^{n+k}$. It is very interesting to note that many aspects of closed minimal submanifolds in the unit sphere have an analogous one for free boundary minimal submanifolds in the unit ball. In this sense, the Ambrozio-Nunes' Theorem can be seen as the analogous result of \cite{CdCK}, \cite{Lawson}, and \cite{Simons}. 

Our goal in this paper is to obtain a topological gap theorem for compact free boundary (not necessarily minimal) submanifolds in the unit ball. In order to state our theorems, we recall that the traceless second fundamental form is defined as $\Phi(u,v)=A(u,v)-\langle u,v\rangle\vec{H}$, for $u,v\in T_x M$, and  $x\in M$, where $\vec{H}$ is the mean curvature vector of $M$. In particular, $\|\Phi\|^2 = \|A\|^2- n\|\vec H\|^2$.

There exists a significant difference between the $2$-dimensional case and the $n$-dimensional ones, for $n\geq 3$, in terms of the technique employed to obtain the results. In fact, the $2$-dimensional case follows from well-known properties of free boundary surfaces, while for higher dimensions, we need to apply more sophisticated tools. We have the following theorems.

\begin{theorem}\label{therem1.2}
Let $\Sigma^2$ be a free boundary compact orientable surface immersed in $\mathbb{B}^{2+k}$, for any positive integer $k$. If $\|\Phi\|^2\le 2$, then $\Sigma^2$ is topologically a disk.
\end{theorem}

If $\Sigma^2$ is minimal, we can improve the constant of Theorem \ref{therem1.2} and, in virtue of a result due to Fraser and Schoen \cite{FS}, obtain the following rigidity result.

\begin{corollary}\label{corollary.1.3} 
Let $\Sigma^2$ be a free boundary compact orientable surface minimally immersed in $\mathbb B^{2+k}$, for any positive integer $k$. If $\|A\|^2\le 4$, then $\Sigma^2$ is the flat equatorial disk.
\end{corollary}

In higher dimensions, we prove that there are no nontrivial harmonic $p$-forms on $M^n$ with either Neumann or Dirichlet condition on the boundary. So, using the Hodge-de Rham Theorem, we actually have the following result.

\begin{theorem}\label{theorem.main}
Let $M^n$ be a compact oriented submanifold immersed in $\B^{n+k}$, with $n\ge3$, which is free boundary and has flat normal bundle. If $\|\Phi\|^2<\frac{np}{n-p}$, for some positive integer $p\le\left\lfloor\frac{n}{2}\right\rfloor$, then the $p$th and the $(n-p)$th cohomology groups of $M^n$ with real coefficients vanish, that is, $H^p(M;\mathbb R)=H^{n-p}(M;\mathbb R)=0.$ In particular, if $\|\Phi\|^2<\frac{n}{n-1}$, then all cohomology groups $H^q(M;\mathbb{R})$, with $q=1,\ldots,n-1$, vanish and $M$ has only one boundary component.
\end{theorem}

In order to prove this theorem, we apply the well-established B\"ochner's technique together with an appropriate estimate for the Weitzenb\"ock tensor of $M^n$ in terms of its extrinsic geometry. To deal with the boundary term, we use a Hardy type inequality for submanifolds which was recently discovered by Batista, Mirandola and the third author \cite{BatistaMirandolaVitorio}. 

When $M^n$ is minimal, we can improve the constant obtained in Theorem \ref{theorem.main}.

\begin{theorem}\label{theorem.main2}
Let $M^n$ be a compact oriented submanifold minimally immersed in $\B^{n+k}$, with $n\ge3$, which is free boundary and has flat normal bundle. Given a positive integer $p\le\left\lfloor\frac{n}{2}\right\rfloor$, we have the following assertions:
\begin{enumerate}
\item If $\|A\|^2<\frac{n^2}{2(n-p)}$, then $H^p(M;\R)$ vanishes. If, additionally, $p=\left\lfloor\frac{n}{2}\right\rfloor$, then $H^{n-p}(M;\R)$ also vanishes.
\item If $\|A\|^2\le\frac{(n-p+1)n^3}{4p(n-p)^2}$ and $1\le p\le\left\lfloor\frac{n}{2}\right\rfloor-1$, then $H^{n-p}(M;\R)$ vanishes.
\end{enumerate}
In particular, if $\|A\|^2<\frac{n^2}{2(n-1)}$, then all cohomology groups $H^q(M;\mathbb{R})$, with $q=1,\ldots,n-1$, vanish and $M$ has only one boundary component.
\end{theorem}

Our theorems lead us to the following questions:\\

{\bf Open questions:} {\it Do Theorems \ref{theorem.main} and \ref{theorem.main2} hold
without the condition on the flatness of the normal bundle? What are the best constants in such
theorems?}\\



\bigskip
{\bf Acknowledgments:} 
The authors are grateful to  Levi Lima and  Ezequiel Barbosa 
 for their interest and helpful discussions about this work.  
The authors were partially supported by CNPq-Brazil, CAPES-Brazil and FAPEAL-Brazil.

\section{Preliminaries}

Let $M^n$ be a compact Riemannian $n$-manifold with nonempty boundary. Let denote by $\Omega^p(M)$ the space of differential $p$-forms on $M$, $d:\Omega^p(M)\to\Omega^{p+1}(M)$ the exterior derivative, and $d^*:\Omega^p(M)\to\Omega^{p-1}(M)$ the codifferential, 
which can be written in terms of the Hodge star operator on $M$ as $d^*=(-1)^{n(p+1)+1}*d*$. We say that $\omega\in\Omega^p(M)$ is {\em harmonic} if $d\omega=0$ and $d^*\omega=0$ on $M$, that is, $\omega$ is closed and coclosed. 
A harmonic $p$-form $\omega$ on $M$ is called {\em tangential} if 
$$
i_\nu\omega=0\mbox{ on }\p M
$$
and {\em normal} if 
$$
\nu\wedge\omega=0\mbox{ on }\p M.
$$
We can consider the following subspaces of $\Omega^p(M)$:

\begin{eqnarray*}
\Har_N^p(M)=\{\omega\in\Omega^p(M);\omega\mbox{ is harmonic and tangential}\}
\end{eqnarray*}
and
\begin{eqnarray*}
\Har_T^p(M)=\{\omega\in\Omega^p(M);\omega\mbox{ is harmonic and normal}\}.
\end{eqnarray*}

It is well-known that the Hodge star operator gives an isomorphism between $\Har_T^{n-p}(M)$ and $\Har_N^p(M)$. Then, using the Hodge-de Rham Theorem (see \cite[Theorem 3]{ACS}), we have
$$
\Har_T^{n-p}(M)\simeq\Har_N^{p}(M)\simeq H^p(M;\mathbb R).
$$
An important fact about $\Har_T^1(M)$ is that $\dim\Har_T^1(M)\geq r-1$, where $r$ is the number of boundary components of $M$ (see \cite[Lemma 4]{ACS}).

Now, let us present some tools which are going to be of use. We start by recalling the integral version of the Weitzenb\"ock formula for manifolds with umbilical boundary (see, for instance, \cite{Levi} and \cite{Yano}). Note that this is exactly the case when $M$ is a free boundary submanifold in the Euclidean unit ball.

\begin{lemma}[Weitzenb\"ock formula]
If $\p M$ is totally umbilical in $M^n$ with second fundamental form $B=I$, then
\begin{eqnarray*}
&\ds\int_M\|\nabla\omega\|^2+\langle\mathcal R_p(\omega),\omega\rangle=-\alpha\int_{\p M}\|\omega\|^2,&
\end{eqnarray*}
where $\alpha=p$ or $\alpha=n-p$, depending whether $\omega\in\Har_N^p(M)$ or $\omega\in\Har_T^p(M)$, respectively. Here, $\mathcal R_p$ represents the Weitzenb\"ock tensor acting on $p$-forms.
\end{lemma}

Another useful result is the refined Kato's inequality for harmonic forms (see, for instance, \cite{CalderbankGauduchonHerzlich} and \cite{Herzlich}). 
\begin{lemma}[Refined Kato's inequality]
If $\omega$ is a harmonic $p$-form on $M^n$, with $1\le p\le \left\lfloor\frac{n}{2}\right\rfloor$, then
\begin{eqnarray*}
&\ds\|\nabla\omega\|^2\ge\frac{n-p+1}{n-p}\|\nabla\|\omega\|\|^2.&
\end{eqnarray*}
\end{lemma}

Now, we are going to present two results in the context of submanifolds. Let $\B^{n+k}$ be the closed unit ball in $\R^{n+k}$ centered at the origin. Consider a compact oriented immersed submanifold $M^n$ in $\B^{n+k}$ with nonempty boundary $\p M$ and denote by $X$ the unit vector normal to $\p\B^{n+k}=\mathbb{S}^{n+k-1}$ which is outward pointing. Denote by $\nu$ the conormal vector field along $\partial M$, that is, the unit vector normal to $\p M$ and tangent to $M$  which points to the outside of $M$. In this setting, saying that $M^n$ is {free boundary} in $\B^{n+k}$ is equivalent to say that $\nu=X$ along $\p M$. 

Denote by $A$ the second fundamental form of $M$ in $\B^{n+k}$, by $\vec{H}$ the mean curvature vector of $M$ with respect to $A$ and by $\Phi$ the traceless part of $A$, i.e., 
$$
\Phi(u,v)=A(u,v)-\langle u,v\rangle\vec{H},\,\,u,v\in T_x M,\,\, x\in M.
$$

If $M^n$ is a compact oriented immersed submanifold in $\B^{n+k}$ which is free 
boundary and has dimension $n\ge3$, then it holds a Hardy type inequality on $M^n$.
In fact, from a result due to Batista, Mirandola and the third author (see \cite[Theorem 3.2]{BatistaMirandolaVitorio}), we have 
\begin{eqnarray*}
\frac{(n-\gamma)^p}{p^p}\int_M\frac{u^p}{r^\gamma}+\frac{\gamma(n-\gamma)^{p-1}}{p^{p-1}}\int_M\frac{u^p}{r^\gamma}\|\bar\nabla r^\perp\|^2\ \ \ \ \ \ \ \ \ \ \ \ \ \ \ \ \ \ \\ 
\ \ \ \ \ \ \ \ \ \ \ \ \ \ \ \ \ \ \le\frac{1}{p^p}\int_M\frac{1}{r^{\gamma-p}}\|p\nabla u+nu\vec{H}\|^p
+\frac{(n-\gamma)^{p-1}}{p^{p-1}}\int_{\p M}\frac{u^p}{r^{\gamma-1}}
\end{eqnarray*}
for all nonnegative function $u\in C^1(M)$, $p\in[1,\infty)$, and $\gamma\in(-\infty,n)$, where $r$ is the distance function to the origin in $\R^{n+k}$, $\bar\nabla r$ denotes the gradient of $r$ in $\R^{n+k}$, and $\nabla u$ is the gradient of $u$ in $M^n$.

Taking $p=2$, $\gamma=0$ and observing that $r\le1$ on $M^n\subset\B^{n+k}$ 
and $r=1$ on $\p M\subset\mathbb{S}^n$, we have the following

\begin{lemma}[Batista-Mirandola-Vit\'orio]\label{ineq}
If $M^n$ is a compact oriented immersed  submanifold in $\B^{n+k}$ 
which is free boundary and has dimension $n\ge3$, then
\begin{eqnarray*}
&\ds\int_M u^2\le\left(\frac{2}{n}\right)^2\int_M\|\nabla u\|^2+\int_Mu^2\|\vec{H}\|^2+\frac{2}{n}\int_{\p M}u^2&
\end{eqnarray*}
for all nonnegative function $u\in C^1(M)$.
\end{lemma}

Finally, we are going to use an extrinsic estimate for the curvature term which appears in the Weitzenb\"ock formula (see \cite{HeziLin}).

\begin{lemma}[Lin]\label{lemma.x}
If $M^n$ is immersed in $\R^{n+k}$ with flat normal bundle, then
\begin{eqnarray*}
&\!\!\!\!\ds\langle\mathcal R_p(\omega),\omega\rangle\ge\!\left(\!p(n-p)\|\vec{H}\|^2-\frac{p(n-p)}{n}\|\Phi\|^2-|n-2p|\sqrt{\frac{p(n-p)}{n}}\|\vec{H}\|\|\Phi\|\!\right)\!\|\omega\|^2.&
\end{eqnarray*}
\end{lemma}

\section{Proof of Theorem \ref{theorem.main}}

Fix $1\le p\le\left\lfloor\frac{n}{2}\right\rfloor$ and $\omega\in\Har_N^p(M)$ or $\omega\in\Har_T^p(M)$, and define $u=\|\omega\|$. It follows from the Weitzenböck formula and the refined Kato's inequality that 
\begin{eqnarray*}
-\alpha\int_{\p M}u^2&=&\int_M\|\nabla\omega\|^2+\int_M\langle\mathcal{R}_p(\omega),\omega\rangle\\
&\ge&\frac{n-p+1}{n-p}\int_M\|\nabla u\|^2+\int_M\langle\mathcal{R}_p(\omega),\omega\rangle,
\end{eqnarray*}
where $\alpha=p$ or $\alpha=n-p$, 
depending whether $\omega\in\Har_N^p(M)$ or $\omega\in\Har_T^p(M)$, respectively. 
Using Lemma \ref{lemma.x}, we obtain
\begin{eqnarray}\label{des}
\nonumber-\alpha\int_{\p M}u^2&\ge&\frac{n-p+1}{n-p}\int_M\|\nabla u\|^2+p(n-p)\int_Mu^2\|\vec{H}\|^2\\ 
&&-\frac{p(n-p)}{n}\int_Mu^2\|\Phi\|^2-(n-2p)\sqrt{\frac{p(n-p)}{n}}\int_Mu^2\|\vec{H}\|\|\Phi\|.
\end{eqnarray}
Fix $\varepsilon>0$ and observe that $$\|\vec{H}\|\|\Phi\|\le\frac{\varepsilon}{2}\|\vec{H}\|^2+\frac{1}{2\varepsilon}\|\Phi\|^2.$$
Then,
\begin{eqnarray*}
-\alpha\int_{\p M}u^2&\ge&\frac{n-p+1}{n-p}\int_M\|\nabla u\|^2+A\int_Mu^2\|\vec{H}\|^2-B\varphi^2\int_Mu^2,
\end{eqnarray*}
where
\begin{eqnarray*}
&\ds A=A(n,p,\varepsilon)=p(n-p)-\frac{(n-2p)\varepsilon}{2}\sqrt{\frac{p(n-p)}{n}},&\\
&\ds B=B(n,p,\varepsilon)=\frac{p(n-p)}{n}+\frac{(n-2p)}{2\varepsilon}\sqrt{\frac{p(n-p)}{n}},&
\end{eqnarray*}
and $\varphi=\sup_M\|\Phi\|$. Therefore, using Lemma \ref{ineq}, we have
\begin{eqnarray*}
0&\ge&\left(\alpha-\frac{2B}{n}\varphi^2\right)\int_{\p M}u^2+\left(\frac{n-p+1}{n-p}-\frac{4B}{n^2}\varphi^2\right)\int_M\|\nabla u\|^2\\
&&+\left(A-B\varphi^2\right)\int_Mu^2\|\vec{H}\|^2.
\end{eqnarray*} 
It follows from the above inequality that if 
\begin{itemize}
\item $\ds\varphi^2<\frac{n\alpha}{2B}=\phi_1$,
\item $\ds\varphi^2\le\frac{(n-p+1)n^2}{4(n-p)B}=\phi_2$, and
\item $\ds\varphi^2\le\frac{A}{B}=\phi_3$,
\end{itemize}
then $u=\|\omega\|=0$ on $\p M$, which implies that $\omega=0$ since $\omega$ is harmonic. 

Now, to finish the proof of Theorem \ref{theorem.main}, we are going to prove that $$\frac{np}{n-p}=\phi_3(\varepsilon_1)=\max_{\varepsilon>0}\phi_3\le\min\{\phi_1(\varepsilon_1),\phi_2(\varepsilon_1)\},$$ for $$\varepsilon_1=\sqrt{\frac{np}{n-p}}.$$ 

\begin{claim}\label{claim.3} 
$\ds\frac{np}{n-p}=\phi_3(\varepsilon_1)=\max\limits_{\varepsilon>0}\phi_3$.
\end{claim}

First, suppose that $2p<n$. It is not difficult to see that there exists a unique $a=a(n,p)$ such that $A(n,p,a)=0$. In fact, $$a=\frac{2\sqrt{np(n-p)}}{n-2p}.$$ Furthermore, $A(n,p,\varepsilon)>0$ for $\varepsilon\in(0,a)$ and $A(n,p,\varepsilon)<0$ for $\varepsilon>a$. Also, $\lim\limits_{\varepsilon\to 0^+}\phi_3=0$. Then, to calculate $$\max_{\varepsilon>0}\phi_3=\max_{\varepsilon\in(0,a)}\phi_3,$$ it is sufficient to find the critical points of $\varepsilon\longmapsto\phi_3(n,p,\varepsilon)$ on the interval $(0,a)$. Define $$b=\frac{(n-2p)}{2}\sqrt{\frac{p(n-p)}{n}}$$ and $$c=\frac{p(n-p)}{n},$$ and observe that $$\phi_3=\frac{A}{B}=\frac{nc\varepsilon-b\varepsilon^2}{c\varepsilon+b}.$$ A straightforward calculation gives that the unique critical point of $\phi_3$ on the interval $(0,+\infty)$ is given by $$\varepsilon_1=\frac{-b+\sqrt{b^2+nc^2}}{c}.$$ Observing that  
\begin{eqnarray*}
b^2+nc^2&=&\frac{p(n-p)(n-2p)^2}{4n}+\frac{p^2(n-p)^2}{n}\\
&=&\frac{pn(n-p)}{4},
\end{eqnarray*}
we have
\begin{eqnarray*}
\varepsilon_1&=&\frac{-b+\sqrt{b^2+nc^2}}{c}\\
&=&\frac{n}{p(n-p)}\left(-\frac{(n-2p)}{2}\sqrt{\frac{p(n-p)}{n}}+\frac{\sqrt{pn(n-p)}}{2}\right)\\
&=&\frac{n}{2\sqrt{p(n-p)}}\left(-\frac{n-2p}{\sqrt{n}}+\sqrt{n}\right)\\
&=&\sqrt{\frac{np}{n-p}}.
\end{eqnarray*}
Then, evaluating $A$ and $B$ at $\varepsilon=\varepsilon_1$, we obtain
\begin{eqnarray*}
A(n,p,\varepsilon_1)=\frac{np}{2}
\end{eqnarray*}
and 
\begin{eqnarray*}
B(n,p,\varepsilon_1)=\frac{n-p}{2}.
\end{eqnarray*}
Thus,
$$\max_{\varepsilon>0}\phi_3=\phi_3(\varepsilon_1)=\frac{np}{n-p}.$$

If $2p=n$, then $A=p(n-p)=p^2$ and $B=\frac{p(n-p)}{n}=\frac{p}{2}$. Therefore, $\phi_3$ is constant equal to 
\begin{eqnarray*}
\frac{A}{B}=2p=\frac{np}{n-p}.
\end{eqnarray*}

\begin{claim}
$\phi_1(\varepsilon_1)\ge\phi_3(\varepsilon_1)$.
\end{claim}

Observe that $2B(n,p,\varepsilon_1)=n-p$. Then, $\phi_1(\varepsilon_1)\ge\phi_3(\varepsilon_1)$ is equivalent to $\alpha\ge p$, which is clearly true because we are assuming that $p\le\left\lfloor\frac{n}{2}\right\rfloor\le\frac{n}{2}$.

\begin{claim}
$\phi_2(\varepsilon_1)>\phi_3(\varepsilon_1)$.
\end{claim}
Using that $2B(n,p,\varepsilon_1)=n-p$ and $p\le\frac{n}{2}$, we obtain
\begin{eqnarray*}
\phi_2(\varepsilon_1)=\frac{(n-p+1)n^2}{2(n-p)^2}>\frac{n^2}{2(n-p)}\ge\frac{np}{n-p}=\phi_3(\varepsilon_1).
\end{eqnarray*}

\section{Proof of Theorem \ref{theorem.main2}}

We start from inequality (\ref{des}) assuming that $\vec H=0$, that is,
$$
-\alpha\int_{\p M}u^2\ge\frac{n-p+1}{n-p}\int_M\|\nabla u\|^2
-\frac{p(n-p)}{n}\int_Mu^2\|\Phi\|^2.
$$
Using Lemma \ref{ineq} in this case, we get
\begin{eqnarray*}
0\geq\left(\alpha-\frac{2p(n-p)}{n^2}\varphi^2 \right)\int_{\partial M}u^2 
+\left(\frac{n-p+1}{n-p} - \frac{4p(n-p)}{n^3}\varphi^2 \right)\int_{ M}\|\nabla u\|^2,
\end{eqnarray*}
and again, the theorem follows if 
\begin{itemize}
\item $\ds\varphi^2<\frac{\alpha n^2}{2p(n-p)}=\phi_4$, and
\item $\ds\varphi^2\le\frac{(n-p+1)n^3}{4p(n-p)^2}=\phi_5$.
\end{itemize}

Studying the solutions of the algebraic equation $\phi_4=\phi_5$, we obtain:

\begin{enumerate}
\item If $\alpha=p$, then $\phi_4<\phi_5$, and it corresponds to the first part of assertion (a). 

\item If $\alpha=n-p$ and $p=\left\lfloor\frac{n}{2}\right\rfloor$, then $\phi_4<\phi_5$, and it corresponds to the second part of assertion (a).

\item In the remaining cases we have $\phi_5<\phi_4$, and they correspond to assertion (b).

\end{enumerate}

\section{Proof of Theorem \ref{therem1.2} and Corollary \ref{corollary.1.3}}

Let $\Sigma$ be a free boundary compact orientable surface immersed in the Euclidean unit ball $\mathbb{B}^{2+k}$. If $K$ represents the Gaussian curvature of $\Sigma$ and $k_g$ represents the geodesic curvature of $\partial\Sigma$ in $\Sigma$, the Gauss-Bonnet Theorem says that 
\begin{eqnarray*}
\int_{\Sigma} K+\int_{\partial\Sigma}k_g=2\pi\chi(\Sigma)=2\pi(2-2g-r),
\end{eqnarray*}
where $\chi(\Sigma)$ is the Euler characteristic of $\Sigma$, $g$ is the genus of $\Sigma$, and $r$ is the number of connect components of $\partial\Sigma$. From the assumption that $\Sigma$ is free boundary in $\mathbb{B}^{2+k}$, we have
\begin{itemize}
\item $k_g\equiv 1$;
\item $|\partial\Sigma|=2\int_\Sigma(1+\langle\vec{H},x\rangle)$,
\end{itemize}
where $|\partial\Sigma|$ represents the length of $\partial\Sigma$. On the other hand, the Gauss equation yields $K=\|\vec{H}\|^2-\frac{1}{2}\|\Phi\|^2$. Thus, using the properties above and the Gauss-Bonnet Theorem, we obtain
\begin{equation} \label{gbo}
\int_\Sigma\left(2+2\langle\vec{H},x\rangle+\|\vec{H}\|^2-\frac{1}{2}\|\Phi\|^2\right)=2\pi(2-2g-r),
\end{equation}
which can be rewritten in the form
\begin{equation}\label{gbfc}
\int_\Sigma\left(1-\|x\|^2+\|x+\vec{H}\|^2+1-\frac{1}{2}\|\Phi\|^2\right)=2\pi(2-2g-r).
\end{equation}

Note that, if $\|\Phi\|^2\leq 2$, then the left hand side of \eqref{gbfc} is positive, which implies that $g=0$ and $r=1$. This concludes the proof of Theorem \ref{therem1.2}. 

In the minimal case, equation \eqref{gbo} simplifies to
$$
\int_\Sigma  \left( 2 -\frac{1}{2}\|A\|^2\right)=2\pi (2-2g-r).
$$
Therefore, if $\|A\|^2\le 4$, then $2-2g-r\ge 0$. So, we have two possible cases:
\begin{enumerate}
\item[\rm (i)] $2-2g-r=1$. In this case, $\Sigma$ is topologically a disk and the result follows from the Fraser-Schoen's Theorem \cite{FS}.
\item[\rm (ii)] $2-2g-r=0$. In this case, $\|A\|^2\equiv 4$, and $K\equiv-2$, which contradicts Theorem 6 in \cite{Yau1974}.
\end{enumerate}

\bibliographystyle{amsplain}
\bibliography{bibliography}

\end{document}